\documentclass[a4paper,11pt]{article}
\usepackage{soul}
\usepackage{amsmath, amsthm, amssymb}
\usepackage[T1]{fontenc}
\usepackage[latin1]{inputenc}

\newtheorem{theo}{Theorem}[section]

\newtheorem{lemma}[theo]{Lemma}
\newtheorem{examp}[theo]{Example}
\newtheorem{examps}[theo]{Examples}

\newtheorem{rmrk}[theo]{Remark}
\newtheorem{proposition}[theo]{Proposition}

\newtheorem{cor}[theo]{Corollary}

\newenvironment{remark}{\begin{rmrk}\rm}{\end{rmrk}}
\newenvironment{example}{\begin{examp}\rm}{\end{examp}}
\newenvironment{examples}{\begin{examps}\rm}{\end{examps}}

\newcommand{\N}{{\mathbb N}}

\newcommand{\pal}{{\rm Pal}}
\newcommand{\ipal}{{\rm IPal}}

\newcommand{\bw}{{\bf w}}
\newcommand{\bt}{{\bf t}}
\newcommand{\bs}{{\bf s}}

\newcommand{\A}{{\mathcal A}}

\newcommand{\V}{{\mathcal V}}
\newcommand{\id}{\rm id}
\newcommand{\h}{\mathcal H}
\newcommand{\e}{{\rm E}}

\begin{document}


\title{On the fixed points of the iterated pseudopalindromic closure}

\author{D. Jamet \footnotemark[1]
\and
G. Paquin\footnotemark[4]   \footnotemark[2]
\and 
G. Richomme \footnotemark[3]
\and
L. Vuillon \footnotemark[2]
}

\footnotetext[1]{LORIA - Universit\'e Nancy 1 - CNRS, Campus Scientifique, BP 239, 54506 Vandoeuvre-l\`es-Nancy, France, Damien.Jamet@loria.fr}
\footnotetext[2]{Laboratoire de math\'ematiques, CNRS UMR 5127, Universit\'e de Savoie, 73376 Le Bourget-du-lac cedex, France, [Genevieve.Paquin, Laurent.Vuillon]@univ-savoie.fr }
\footnotetext[3] {UPJV, Laboratoire MIS,
33, Rue Saint Leu, 80039 Amiens cedex 0, France,
Gwenael.Richomme@u-picardie.fr}
\footnotetext[4]{with the support of FQRNT (Qu\'ebec)}
\date{\today}

\maketitle

\begin{abstract} 
First introduced in the study of the Sturmian words, the iterated palindromic closure was recently generalized to pseudopalindromes. This operator allows one to construct words with an infinity of pseudopalindromic prefixes, called pseudostandard words. We provide here several combinatorial properties of the fixed points under the iterated
pseudopalindromic closure. 
\end{abstract}


\section{Introduction}

The Sturmian words form a well-known class of infinite words over a $2$-letter alphabet that occur in many different fields, for instance in astronomy,  symbolic dynamics, number theory, discrete geometry, cristallography, and of course, in combinatorics on words (see \cite{Lothaire2002} for instance).  Depending on the context of the study, these words have many equivalent characterizations. In discrete geometry, they are exactly the words that approximate a discrete line having irrational slope, using horizontal and diagonal moves. In symbolic dynamics, they are obtained by the exchange of $2$ intervals. They are also known as the balanced aperiodic infinite words over a $2$-letter alphabet. A subclass of the Sturmian words is formed by the standard Sturmian ones.  For each Sturmian word, there exists a standard one having the same language, i.e. the same set of factors. A standard Sturmian word is, in a sense, the representative of all Sturmian words having the same language. All the words in this subclass can be easily obtained by a construction called the  iterated palindromic closure \cite{deLuca1997}. This operation is a bijection between standard Sturmian words and infinite words over a $2$-letter alphabet that do not end by the repetition of a unique letter. 

On the other side, some fixed points of function are famous in combinatorics on words. As an example, the  self-generating word introduced in \cite{Kolakoski1965}, called the Kolakoski word,  is the fixed point under the run-length encoding function and raised some challenging problems. For instance, we still do not know what are its letter frequencies, if they exist. The recurrence of the Kolakoski word as well as the closure of its set of factors under complementation or reversal are other open problems. 

In this context, it is a natural problem to try to characterize the fixed points under the iterated palindromic closure operator, and more generally, under the iterated pseudopalindromic closure operator, recently introduced in \cite{deLucaDeLuca2006}. In this paper, we study these words and show some of their  properties. It is organized as follow. In Section $2$, we first  recall what is the iterated palindromic closure operator and then, in Section $3$, we introduced the iterated pseudopalindromic closure operator, which generalized the first one using a generalization of a palindrome. 
In Section $4$, we prove  the existence of fixed points under the iterated pseudopalindromic closure operator and we show them explicitly: there are $3$ families of fixed points.  
Finally, in Section $5$, we give some of their combinatorial properties.

Let us note that by lack of space, many proofs are skipped. Nevertheless we provide the main intermediate steps. Note also that we assume the reader is familiar with the notions and notations of Combinatorics on Words (see \cite{Lothaire2002} for instance). In all the paper $\A$ denotes an alphabet.

\section{Iterated palindromic closure}
Sturmian words may be defined in many equivalent ways (see Chapter 2 in \cite{Lothaire2002} for more details). For instance, they are the non ultimately periodic infinite words over a $2$-letter alphabet  that have the minimal complexity, that is the number of distinct factors of  length $n$  is $(n+1)$. They are also the set of non ultimately periodic binary balanced words. Recall that a binary word $w$ is {\it balanced} if for all factors $f, f'$ having same length, and for all letter $a \in \A$, one has $\left | |f|_a-|f'|_a\right | \leq 1$.

The Sturmian words also are the infinite non ultimately binary words that describe a discrete line.
Recall that the slope of the word $\bf s$ is $\alpha=\lim_{n\rightarrow \infty }|s[1..n]|_b/n$. 

All Sturmian words considered in this paper belong to the particular class of standard Sturmian words for which we now recall the construction using  the iterated palindromic closure operator.
Given a finite word $w$, let us denote by $\pal(w)$ the word obtained iterating the palindromic closure: $\pal(\varepsilon) = \varepsilon$ and $\pal(wa) = (\pal(w)a)^{(+)}$, for all words $w$ and letters $a$.

By the definition of the iterated palindromic closure $\pal$, for any finite word $w$ and letter $a$, $\pal(w)$ is a prefix of $\pal(wa)$. One can then define for any infinite word $\bw = (a[n])_{n \geq 1}$, the infinite word $\ipal(\bw) = \lim_{n \to \infty} \pal(a[1]\cdots a[n]).$ We then say that the word $\bf w$ {\it directs} the word $\ipal(\bf w)$.  From the works of \cite{deLuca1997}, we know that 
$\ipal$ is a bijection between the set of binary infinite words not of the form $ua^\omega$, with $u\in \A^*$ and $a \in \A$, and the set of standard Sturmian words. The word $\bw$ is then  called the {\it directive word}  of the standard Sturmian word $\ipal(\bw)$. Note that words of the form $\ipal(ua^\omega)$ (with $u\in \A^*$ and $a \in \A$) are periodic (see  Lemma~\ref{L:periodic episturmian}  below recalled from \cite{DroubayJustinPirillo2001}).

The $\ipal$ operator is also well defined over a $k$-letter alphabet, with $k\geq 3$. In this case, it is known \cite{DroubayJustinPirillo2001} that $\ipal(\A^\omega)$ is the set of {\it standard episturmian words}, a generalization over a $k$-letter alphabet, $k\geq 3$, of the family of standard Sturmian words (for more details, see \cite{GlenJustin2009TIA}). When $\bw$ is a word over $\A$ containing infinitely often each letter, then $\ipal(\bw)$ is called a \textit{strict} standard episturmian word. The set of strict standard episturmian words corresponds to the set of Arnoux-Rauzy words \cite{ArnouxRauzy1991}.

%

\begin{example} The infinite word $abcabaac\cdots$ directs the standard episturmian word $\bw=\pal(abcabaac\cdots)=\underline a \underline ba \underline caba \underline abacaba \underline b acabaabacaba \underline a\cdots .$
\end{example}

As we will do in the sequel, in the previous example we have underlined the letters in the standard episturmian word corresponding to the letter of its directive word, for the sake of clarity.


\section{Iterated pseudopalindromic closure}

Recently, de Luca and De Luca \cite{deLucaDeLuca2006} have extended the notion of palindrome to what they call {\it pseudopalindrome}, using involutory antimorphism. Let
recall that a map $\V:\A^* \rightarrow \A^*$ is called an {\it antimorphism} of $\A^*$ if for all $u,v \in \A^*$ one has $\V(uv)=\V(v)\V(u)$. Moreover, an antimorphism is {\it involutory} if $\V^2=\id$.
 A trivial involutory antimorphism is the {\it reversal} $\, \, \widetilde {\,}\, $ which we will denote in the sequel by the function $R: \A^* \rightarrow \A^*$, $R(w)=\widetilde{w}$. Any involutory antimorphism $\V$ of $\A^*$ can be constructed as $\V=\tau \circ R= R\circ \tau$, with $\tau$ an involutory permutation of the alphabet $\A$. From now on, in order to describe an involutory antimorphism $\V$, we will then only give the involutory permutation $\tau$ of the alphabet $\A$. The two antimorphisms $E$ and $\h$ defined respectively over $\{a, b\}$ and $\{a, b, c\}$ by\\
\centerline{$E = R \circ \tau$ with $\tau(a) = b$, $\tau(b) = a$,}\\
\centerline{$\h = R \circ \tau$ with $\tau(a) = a$, $\tau(b) = c$, $\tau(c) = b$}
will play, in addition to $R$, an important role in our study. The antimorphism $E$ will be called, as usually, the \textit{exchange antimorphism}. We propose to name antimorphism $\h$ the \textit{hybrid antimorphism}, hence the notation, since it contains both an identity part and an exchange part.


We can now define the generalization of palindromes given in \cite{deLucaDeLuca2006}:  a word $w \in \A^*$ is called a {\it $\V$-palindrome} if it is the fixed point of an involutory antimorphism $\V$ of the free monoid $\A^*$: $\V(w)=w$. When the antimorphism $\V$ is not mentioned, we call it a {\it pseudopalindrome}. 



Similarly to the palindromic closure $^{(+)}$, the {\it $\V$-palindromic closure} of the finite word $u$, also called the   {\it pseudopalindromic closure} when the antimorphism is not specified,  is defined by $u^{\oplus}=sq\V(s)$, where $u=sq$, with $q$ the longest $\V$-palindromic suffix of $u$. The pseudopalindromic closure of $u$ is the shortest pseudopalindrome having $u$ as prefix.

\begin{example}Over the alphabet $\{a,b\}$, since the longest $E$-palindromic suffix of $w = aaba$ is $ba$, $w^\oplus=aaba\cdot E(aa)=aababb$. 
\end{example}


Extending the $\pal$ operator to pseudopalindrome, the $\pal_ \V$ operator is naturally defined by $\pal_ \V(\varepsilon)=\varepsilon$ and $\pal_\V(wa)=(\pal_\V(w)a)^\oplus$, for $w\in \A^*$ and $a \in \A$. Then, for $\bw\in \A^\omega$, $\ipal_\V(\bw)=\lim_{n \rightarrow \infty} \pal_\V(w[1] \cdots w[n])$. This limit exists since by the definition of $\pal_\V$, for any involutory antimorphism $\V$, $w \in \A^*$ and $a\in \A$, $\pal_\V(w)$ is a  prefix of $\pal_\V(wa)$.
 The infinite word obtained by the $\ipal_\V$ operator is a {\it $\V$-standard word}, also called a {\it pseudostandard word} when the antimorphism is not specified. This new class of words  is a general one that includes the standard Sturmian and the standard episturmian ones and was first introduced in \cite{deLucaDeLuca2006}.


\begin{example} Over $\A=\{a,b,c\}$, 
the $\h$-standard word directed by $(abc)^\omega$ is $\ipal_ \h ((abc)^\omega )= \underline a \underline bca \underline cbabca \underline abcacbabca\underline b cacbabcaabcacbabca \cdots.$
\end{example}


\section{Existence of fixed points}

In this section, we prove the existence of fixed points over the iterated pseudopalindromic closure and we show which forms they can have. We denote  naturally $\ipal_\V^0(\bw) = \bw$ and $\ipal^n_\V(\bw) = \ipal_\V(\ipal_\V^{n-1}(\bw))$, for any $\bw \in \A^\omega$, involutary antimorphism $\V$ and $n  \geq 1$. Let us see some examples of the iteration of the $\ipal_\V$ operator over infinite words.
 
\begin{examples} Over a $2$-letter alphabet $\A=\{a,b\}$, there are only two possible involutory antimorphisms: the reversal antimorphism $R$ and the exchange antimorphism $E$. 
Let us consider for instance the iteration of the $\ipal_R$ operator over the word $\bw=abx\cdots$, with $x\in \A$ (the iteration of $\ipal_\e$ leads to similar remarks):
\vspace{-0.2cm}
  \begin{eqnarray*}
\ipal_R(abx\cdots)&=&\underline a \underline ba \underline x\cdots \\
        \ipal^2_R(abx\cdots)&=&\underline a \underline ba \underline aba \underline x \cdots  \\
        \ipal^3_R(abx\cdots)&=&\underline a \underline ba \underline aba \underline aba \underline baabaaba \underline ababaabaaba \underline x \cdots . 
        \end{eqnarray*}
\end{examples}
\vspace{-0.2cm}
We see that the position of the letter $x$ of the directive word $\bw$ in $\ipal^k_R(\bw)$ grows with the value of $k$. We also observe that the common prefix of $\ipal_R^k(\bw)$ and $\ipal_R^{k+1}(\bw)$ also seems to grow with $k$.  It appears that only a short prefix of $\bw$ is necessary to determine the word obtained by infinitely iterating the $\ipal_R$ operator. Theorem \ref{T:existence} is a direct corollary of this observation and of the following lemmas that can be proved inductively.

\begin{lemma} \label{L:upref} Let $\V=R\circ \tau$ be an involutory antimorphism and let $(u_k)_{k\geq 1}$ be a sequence of finite words defined by 
\begin{center}
$u_1= \left \{ \begin{tabular}{ll}
$a^nb$ &\textnormal{if} $\tau(a)=a$,\\
$a$ &\textnormal {if} $\tau(a)=b$,
\end{tabular} \right. $ \\ 
\end{center}
and for $k\geq 2$, $u_k=\pal_\V(u_{k-1})$, with $a\neq b\in \A$, $n \geq 1$. Then $\lim_{k\rightarrow \infty} u_k$ exists. 
\end{lemma}



\begin{lemma}\label{L:uprefPal} Let $(u_k)_{k\geq 1}$ be the same sequence as in Lemma \ref{L:upref} and let consider an infinite word $\bw$ having $u_1$ as prefix. Then for all $k\geq 1$, $u_k$ is a proper prefix of $\ipal_\V^{k-1}(\bw)$. 
\end{lemma}


\begin{theo}[and definition] \label{T:existence} Over a $k$-letter alphabet, with $k \geq 2$, there are $3$ kinds of fixed points having at least $2$ different letters, only depending on the first letters of the word and the involutory antimorphism $\V = R \circ \tau$ considered. 
\begin{enumerate}
\item When $\tau(a) = a$ and $\tau(b)=b$, with $a\neq b$, for all $n \geq 1$, $\ipal_\V$ has a unique fixed point beginning with $a^nb$, denoted ${\bs}_{R,n,a,b}$, which equals
$${\bs}_{R,n,a,b}=\lim_{i\rightarrow \infty} \pal^i(a^nb)=\underline a^n \underline b a^n (\underline a ba^n )^{n+1}\underline b(a^{n+1}b)^{n+1}a^n\underline a\cdots .$$
\vspace{-0.7cm}
\item When $\tau(a)=a$ and $\tau(b)=c$ for pairwise different letters $a, b, c$, for all $n \geq 1$, $\ipal_\V$ has a unique fixed point beginning with $a^nb$, denoted by ${\bs}_{\h,n,a,b,c}$, which equals
$$\bs_{\h,n,a,b,c}=\lim_{i \rightarrow \infty} \pal_\h^i(a^nb)=\underline a^n\underline bca^n\underline c ba^nbca^n (\underline a bca^ncba^nbca^n)^n\underline c \cdots .$$
\vspace{-0.7cm}
\item When $\tau(a)=b$ and $\tau(b)=a$, with $a\neq b$, $\ipal_\V$ has a fixed point beginning with $a^nb$ only if $n =1$. It is denoted by ${\bs}_{\e,a,b}$ and equals
$$\bs_{\e,a,b}=\lim _{i \rightarrow \infty} \pal_\e^i(a)=\underline ab \underline baab \underline baab \underline abbaabbaab \underline abbaabbaab\underline b\cdots .$$
\end{enumerate}
\end{theo}


Theorem~\ref{T:existence} characterizes all possible fixed points of $\ipal_\V$ except the trivial fixed point of the form $a^\omega$, which is a fixed point for $\ipal_\V$ using any antimorphism $\V=R\circ \tau$ with $\tau(a)=a$. This trivial fixed point corresponds to the words obtained in Theorem \ref{T:existence} {\it 1.} and {\it 2.} with $n=\infty$.

\begin{remark}
Even if there exist many involutory antimorphisms for arbitrary $k$-letter alphabets \cite{deLucaDeLuca2006}, 
fixed points over the $\ipal_\V$ operators contain at most three letters. More precisely, the fixed points over a $3$-letter alphabet $\{a, b, c\}$ starting by $a$ can only be obtained by the antimorphism $\h$ (we recall that $\h=R\circ \tau$, with $\tau(a)=a$ and $\tau(b)=c$). Indeed, $\tau(a)=b$ yields to $\bs_{\e,a,b}$ and $\tau(a)=a$ and $\tau(b)=b$ yields to $\bs_{R,n,a,b}$. Moreover, for the antimorphism $\e$, the fixed point can not start by $a^2$, since $a^2$ is not a prefix of $\pal_\e(a^2)=abab$. 
\end{remark}



\section{Combinatorial properties of the fixed points}

In this section, we consider successively the fixed points ${\bs}_{R,n,a,b}$, $\bs_{\e,a,b}$ and  $\bs_{\h,n,a,b,c}$  of the $\ipal_\V$ operator and we give some of their combinatorial properties. We will see that words  ${\bs}_{R,n,a,b}$ are Sturmian and $\bs_{\e,a,b}$ is related to a Sturmian word, whereas words $\bs_{\h,n,a,b,c}$ cannot be such, since they contain the three letters $a$, $b$ and $c$. This explains why we consider the word $\bs_{\e,a,b}$ before words $\bs_{\h,n,a,b,c}$ contrarily to their order of introduction in Theorem~\ref{T:existence}.

\subsection{Study of the fixed point ${\bs}_{R,n,a,b}$}

Here, we consider the first fixed point of the $\ipal_\V$ operator, with $\V=R$. 
Note that $\ipal_R = \ipal$.
 Before stating our first property, we need the following lemma.

\begin{lemma}[\cite{DroubayJustinPirillo2001}, Theorem 3]\label{L:periodic episturmian}  An infinite word obtained by the $\ipal$ operator is ultimately periodic if and only if its directive word has the form $ua^\omega$, with $u \in \A^*$ and $a\in \A$. 
\end{lemma}

\begin{proposition}\label{p:standard} For a fixed positive $n \in \N$, ${\bs}_{R,n,a,b}$ is not ultimately periodic and consequently, is a standard Sturmian word.
\end{proposition}

\begin{proof} By definition of the word ${\bs}_{R,n,a,b}$,  $(\pal^i(a^nb))_{i \geq 0}$ forms a sequence of prefixes of ${\bs}_{R,n,a,b}$. The sequence of lengthes of these prefixes is strictly increasing by the definition of the Pal operator. Since $ba^n$ is a suffix of $\pal^i(a^nb)$, both letters $a$ and $b$ occur infinitely often in ${\bs}_{R,n,a,b}$. Hence ${\bs}_{R,n,a,b}$ is not of the form $u\alpha^\omega$ for a word $u$ and a letter $\alpha$. Since by its definition, ${\bs}_{R,n,a,b}$ equals its directive word, Lemma~\ref{L:periodic episturmian} implies that ${\bs}_{R,n,a,b}$ is not ultimately periodic.
\end{proof}

Proposition \ref{p:standard} is very useful, since it allows us to use properties of standard Sturmian words in order to characterize the fixed point $\bs_{R,n,a,b}$. Let us recall some of them.

\begin{theo}[\cite{ArnouxRauzy1991}, p. 206]\label{T:Slope Sturmian} Let ${\Delta(\bf w)}=a^{d_1}b^{d_2}a^{d_3}b^{d_4}\cdots$ be the directive word of an infinite standard Sturmian word $\bf w$, with $d_i \geq 1$. Then the slope of  $\bf w$ has the continued fraction expansion $\alpha_\bw=[0;1+d_1,d_2,d_3,d_4, \ldots]$.  
\end{theo}

\begin{theo}{\rm \cite{CrispMoran1993}}\label{T:morphic Sturmian} The standard Sturmian word of slope $\alpha$ is a fixed point of some nontrivial morphism if and only if $\alpha$ is a Sturm number, that is $\alpha$ has a continued fraction expansion of one of the following kinds:
\begin{enumerate}
\item $\alpha=[0;1,a_0,\overline{a_1,\ldots, a_k}]$, with $a_k\geq a_0$,
\vspace{-0.2cm}
\item $\alpha=[0;1+a_0,\overline{a_1,\ldots, a_k}]$, with $a_k\geq a_0\geq 1$.
\end{enumerate}
\end{theo}

%

A wide literature is devoted to the study of these fixed points and the known results about their generating morphism are often used in order to find some of their properties. It is thus natural to wonder if the fixed points of the $\ipal$ operator are also  fixed points of some nontrivial  morphisms. Using Theorems~\ref{T:Slope Sturmian} and \ref{T:morphic Sturmian}, one can see as a quite direct consequence of Proposition~\ref{p:standard}:

\begin{proposition}\label{p:non morphism} For a fixed $n$, ${\bs}_{R,n,a,b}$ is not a fixed point of a nontrivial  morphism.
\end{proposition}


We denote $\alpha_{n,a,b}$ the slope associated  to ${\bf s}_{R,n,a,b}$. One can easily see that the continued fraction expansion $[0; 1+d_1, d_2, \ldots]$ of $\alpha_{n,a,b}$ only contains the letters $0$, $1$, $n$ and $n+1$. Hence:

\begin{lemma}\label{l:bounded} 
The continued fraction expansion of $\alpha_{n,a,b}$ has bounded partial coefficients.
\end{lemma}


Lemma \ref{l:bounded} in itself is not that interesting, but combining it with next lemma allows to get  Proposition \ref{p:repet}. 

\begin{lemma}[\cite{Vandeth2000}, Theorem 17]\label{L:power-free} 
Let $\alpha >0$ be an irrational number with $d_\alpha= [d_0; d_1, d_2, \ldots]$, its continued fraction expansion. Then the standard Sturmian word of slope $\alpha$ denoted ${\bf w}_\alpha$ is $k$-th power-free for some integer $k$ if, and only if, $d_\alpha$ has bounded partial coefficients. Moreover, if $d_\alpha$ has bounded partial coefficients, then ${\bf w}_\alpha$ is $k$-th power-free but not $(k-1)$-th power-free for $k=3+\max _{i\geq0}d_i$.  
\end{lemma}



\begin{proposition} \label{p:repet}
${\bf s}_{R,n,a,b}$ is $(n+4)$-th power-free, but contains $(n+3)$-th powers. 
\end{proposition}

By direct computation, we easily obtain arbitrarily large prefix of the word ${\bf s}_{R,n,a,b}$ for a fixed $n$. The continued fraction expansion of ${\alpha}_{R,n,a,b}$ is then obtained and yields the value of the slope. For $n=1$, we get:
{\small
\begin{eqnarray*}
\alpha_{1,a,b}&=&[0; 2, 1, 2, 1, 2, 1, 1, 1, 2, 1, 2, 1, 2, 1, 1, 1, 2, 1, 2, 1, 2, 1, 1, 1, \ldots ] 
\end{eqnarray*}
}

Whether $\alpha_{n,a,b}$ is transcendental is also an interesting problem and appears as an interesting consequence of Adamczewski and Bugeaud's works \cite{AdamczewskiBugeaud2007}.

\begin{proposition}\label{P:transcendental} For any $n\geq 1$, $\alpha_{n,a,b}$ is transcendental.
\end{proposition}


\begin{theo} {\rm \cite{AdamczewskiBugeaud2007}} \label{t:adam} Let ${\bf a}=(a_\ell)_{\ell \geq 1}$ be a sequence of positive integers. If the word $\bf a$ begins in arbitrarily long palindromes, then the real number $\alpha=[0;a_1,a_2,\ldots, a_\ell, \ldots ]$ is either quadratic irrational or transcendental. 
\end{theo}

\begin{proof} [Proof of Proposition~\ref{P:transcendental}]
 By Proposition \ref{p:standard}, $\bs_{R,n,a,b}$ is not ultimately periodic. Consequently, there are an infinity of occurrences of $a$'s and $b$'s in $\bs_{R,n,a,b}$. Let $$P=\{i \in \N\setminus 0 \, | \,  \bs_{R,n,a,b}[i+1]=a  \} \mbox{ and}$$ $$P'=\{\pal(\bs_{R,n,a,b}[1\ldots i]) \, | \, i \in P\}.$$ Both sets are infinite. Moreover, by its construction, any palindrome in the set $P'$ is followed by an $a$ at its first occurrence in $\bs_{R,n,a,b}$. 
By Theorem~\ref{T:Slope Sturmian} and since $\bs_{R,n,a,b}$ equals its directive word, 
if $a^{i_1}b^{i_2}\cdots b^{i_2}a^{i_1}$ is a palindromic prefix of $\bs_{R,n,a,b}$, 
then the continued fraction expansion of its slope begins by $[0; 1+i_1, i_2, \ldots, i_2, i_1+\rho, \ldots]$, 
for some integer $\rho$.
Moreover by Proposition~\ref{p:standard}, $\bs_{R,n,a,b}$ is standard Sturmian which implies that $\rho = 0$ or $\rho =1$ depending on the next letter occurring in $\bs_{R,n,a,b}$.
By the construction of the palindromes in $P'$, we know that there are all palindromes such that $\rho=1$. That implies that for any $n$, the continued fraction expansion of the slope of $\bs_{R,n,a,b}$ begins by an infinity of palindromes. We conclude using Theorem \ref{t:adam}: since the continuous fraction expansion of the slope is not ultimately periodic, it cannot be quadratic; hence, it is transcendental.
\end{proof}

Notice that the previous proof works since $\bs_{R,n,a,b}$ equals its directive word. Otherwise, the result it not necessarily true. 

\subsection{Study of the fixed point ${\bf  s}_{\e,a,b}$}

We have seen in the previous subsection that since $\bs_{R,n,a,b}$ is a standard Sturmian word, some properties follow directly. Here, we study the fixed point $\bs_{\e,a,b}$.

Recall that Sturmian words are known to be balanced. It is sufficient to consider the letter $a$ and the factors $bb$ and $aa$ to be convinced that $\bs_{\e,a,b}$ is not balanced, and consequently, that it is not a Sturmian word.  

We now recall a powerful result of de Luca and De Luca. For $\V = \tau \circ R$ an involutory antimorphism over an alphabet $\A$, with $\tau$ an involutory permutation of $\A$, $\mu_\V$ is the morphism defined for all $a$ in $\A$, by $\mu(a)=a$ if $a=\tau(a)$ and by $\mu(a)=a\tau(a)$ otherwise.

\begin{theo} \label{t:link} {\rm[\cite{deLucaDeLuca2006}, Theorem 7.1]} For any $\bw \in \A^\omega$ and for any involutory antimorphism $\V$, one has $\ipal_\V(\bw)=\mu_\V(\ipal(\bw)).$
\end{theo}

Since we cannot use the known results about Sturmian words in order to prove combinatorial properties of the fixed point $\bs_{\e,a,b}$, the idea here is to first consider the word $\ipal(\bs_{\e,a,b})$ that will further appear to be standard Sturmian, and then to extend the properties to $\mu_\e(\ipal(\bs_{\e,a,b}))$ which is the fixed point $\bs_{\e,a,b}$, by Theorem \ref{t:link}.

In what follows, $\bw_\e$ will denote $\ipal({\bf  s}_{\e,a,b})$, that is
\begin{equation*}
\bw_\e= \underline a \underline ba \underline ba \underline ababa \underline ababa \underline baababaababa \underline baababaababa \cdots .
\end{equation*}


Notice that here, $\mu_\e$ is the Thue-Morse morphism, that is $\mu_\e(a)=ab$ and $\mu_\e(b)=ba$. 
Note also that $\bs_{\e,a,b}=\mu_\e(\bw_{\e})$, and so that $\bs_{\e,a,b} \in \{ab,ba\}^\omega$. 

Similarly as in the proof of Proposition \ref{p:standard}, one can prove:

\begin{proposition}\label{P:e sturmien} $\bw_\e$ is not ultimately periodic, and consequently, is a Sturmian word.
\end{proposition}

Next lemma claims that $\mu_\V$ morphisms preserve ultimate periodicity. It allows to get Proposition~\ref{p:e not per} which extends the previous proposition to word $\bs_{\e,a,b}$.

\begin{lemma}\label{L:ultim} Let $\V$ be an involutory antimorphism over an alphabet $\A$. 
An infinite word $\bw$ over $\A$ is ultimately periodic if and only if $\mu_\V(\bw)$ is so.
\end{lemma}

\begin{proof}
The "only if" part is immediate.
Assume $\mu_\V(\bw) = uv^\omega$ for words $u \in \A^*$ and $v \in \A^+$. When $v$ begins with a letter $a$ such that $\mu_\V(a) = a$, then $a$ occurs in no word $\mu_\V(b)$ with $b \neq a$, implying that $u = \mu_\V(u')$, $v =\mu_\V(v')$ for some words $u'$, $v'$.
Then $\mu_\V(w) = \mu_\V(u'v'^\omega)$. It is quite immediate that the morphism $\mu_\V$ is injective on infinite words (and also on finite ones). Hence $w = u'v'^\omega$ is ultimately periodic.
Assume now that $v$ begins with a letter $a$ such that $\mu_\V(a) = ab$, we have $\mu_\V(b) = ba$ and neither $a$ nor $b$ occurs in $\mu_\V(c)$ for $c \in \A\setminus\{a, b\}$. Possibly replacing $v$ by $v^2$, we can assume that $|v|_a+|v|_b$ is even. Depending on the parity of $|u|_a+|u|_b$, two cases are possible: $u = \mu_\V(u')$ and $v = \mu_\V(v')$, or, $ua = \mu_\V(u')$ and $a^{-1}va = \mu_\V(v')$. Once again $\mu_\V(w) = \mu_\V(u'v'^\omega)$ and so $w = u'v'^\omega$ is ultimately periodic.
\end{proof}


\begin{proposition}\label{p:e not per}  $\bs_{\e,a,b}$ is not  ultimately periodic.
\end{proposition}

Another way to prove Proposition~\ref{p:e not per} is using the following generalization of Lemma \ref{L:periodic episturmian} to the $\ipal_\V$ operator. 
  
 \begin{proposition}\label{p:pseudo per} 
 Let $\V$ be an involutory antimorphism over an alphabet $\A$.
 An infinite word obtained by the $\ipal_\V$ operator is ultimately periodic if and only if its directive word has the form $u\alpha^\omega$, with $u \in \A^*$ and $\alpha\in \A$.  
\end{proposition}

\begin{proof} 
Assume $\bt = \ipal_\V(\bw)$ is ultimately periodic. By Theorem~\ref{t:link}, $\bt = \mu_\V(\ipal(\bw))$.
Thus Proposition~\ref{p:pseudo per}  appears as a direct corollary of Lemma~\ref{L:ultim} and \ref{L:periodic episturmian}.
\end{proof}

Proposition \ref{p:pseudo per} is interesting by itself, since it generalizes a well-known useful result of Droubay, Justin and Pirillo to pseudostandard words.


By Theorem~\ref{T:Slope Sturmian}, the continued fraction of the slope of $\bw_\e$ is ultimately periodic if and only if its directive word, which is $\bs_{\e,a,b}$ by definition, is ultimately periodic. Hence by Proposition~\ref{p:e not per}, the continued fraction of the slope of $\bw_\e$ is not ultimately periodic which implies next result by Theorem~\ref{T:morphic Sturmian}:


\begin{cor}\label{C:e morphism} $\bw_\e$ is not a fixed point for some nontrivial morphism. 
\end{cor}


A skipped combinatorial proof using Corollary~\ref{C:e morphism} allows to state:

\begin{proposition} ${\bf s}_{\e,a,b}$ is not a fixed point for some non-trivial morphism. \end{proposition}

Let us now consider maximal powers in $\bw_\e$ and $\bs_{\e, a, b}$
\begin{proposition}\label{P:power we} $\bw_\e$ and $\bs_{\e, a, b}$ both contain $4$-th powers, but no $5$-th power words. 
\end{proposition}

\begin{proof} 
By Theorem~\ref{T:Slope Sturmian}, the partial coefficients of the continued fraction of the slope of $\bw_\e$ correspond to the powers of letters on  $\bs_{\e, a, b}$ which can be seen to value $1$ or $2$. 
As a consequence of Lemma~\ref{L:power-free} we deduce the result for $\bw_\e$.


It is known by \cite{Shur2000} that, for all rational $q \geq 2$, a word $\bw$ avoids repetition $u^q$ if and only if $\mu_\e(\bw)$ also avoids them. Thus Proposition~\ref{P:power we} holds for $\bs_{\e,a,b}$.
\end{proof}

Since $\bs_{\e,a,b}$ is not a Sturmian word, it does not have a known geometrical interpretation. 
Thus, the notion of slope does not apply here. However, since $\bs_{\e,a,b} \in \{ab, ba\}^\omega$, we observe that the frequencies of the letters in $\bs_{\e,a,b}$ are both $1/2$.


\subsection{Study of the fixed point ${\bf  s}_{\h,n,a,b,c}$}

Let us now study the properties of the last kind of fixed points. 
Since the words ${\bf  s}_{\h,n,a,b,c}$ do not have a separating letter, they are not episturmian.
As in the previous subsection, let us  denote by $\bw_{\h,n}$ the episturmian word associated by Theorem \ref{t:link} to the fixed point ${\bf s}_{\h, n,a,b,c}$, that is:
\vspace{-0.2cm}
\begin{equation*}
\bw_{\h,n}=\ipal({\bf  s}_{\h,n,a,b,c})= \underline a^n\underline ba^n\underline ca^nba^n\underline a ba^nca^nba^n   \cdots .
\vspace{-0.2cm}
\end{equation*}


As in the proofs of Propositions~\ref{p:standard} and \ref{P:e sturmien}, one can see that the three letters $a, b, c$ occur infinitely often in $\bs_{\h,n,a,b,c}$. Thus by Proposition~\ref{p:pseudo per} and by its construction, $\bf{w}_{\h,n}$ satisfies:

\begin{proposition}\label{P:h sturmien} The words $\bw_{\h,n}$ are not ultimately periodic and are strict  standard episturmian words.
\end{proposition}


Since by definition, $\bs_{\h,n,a,b,c} = \mu_\h(\bw_{\h,n})$, Lemma~\ref{L:ultim} implies:

\begin{proposition}\label{P:S not periodic} $\bs_{\h,n,a,b,c}$ are not ultimately periodic.
\end{proposition}


Let us recall a useful result from Justin and Pirillo.

\begin{proposition}\label{P:JP} {\rm \cite{JustinPirillo2002}} A standard strict episturmian word is a fixed point of a nontrivial morphism if and only if its directive word is periodic.
\end{proposition}

From Propositions \ref{P:S not periodic} and \ref{P:JP}, we get:

\begin{proposition} \label{P:h fixed} $\bw_{\h,n}$ are not fixed points of a nontrivial morphism. 
\end{proposition}





Now we come to repetitions in $\bw_{\h,n}$. 
In \cite{JustinPirillo2002}, Justin and Pirillo provide important tools about fractional powers in episturmian words. We can deduce from their Theorem 5.2 that the critical exponent of any strict  episturmian word $\bs$ having a periodic directive word with the largest block of letters of length  $\ell$, lies between $\ell+2$ and $\ell+3$. In particular $\bs$ is $(\ell+3)$-th power-free but contains an $(\ell+2)$-th power. This can be extended to a larger class of episturmian words, as follows.

%


\begin{proposition}
Let $\bs$ be a strict epistandard word directed by a word $\Delta$ and let $\ell$ denotes the greatest integer $i$ such that $\alpha^i$ is a factor of $\Delta$ with $\alpha$ a letter. Assume $\Delta$ contains at least one factor $aua^\ell va$ with $a$ a letter and $u, v$ non empty words that do not contain the letter $a$. Then $\bs$ is $(\ell+3)$-th power-free but contains an $(\ell+2)$-th power. 
\end{proposition}

\begin{proof}
Let $(v_i)_{i \geq 1}$ be the sequence of prefixes of $\bs$ having a first letter  different from the last letter (it is infinite since $\bs$ is a strict standard episturmian word). For $i \geq 1$, denote $\bs_i$ the standard episturmian word directed by $v_i^\omega$. It is straightforward that $\bs = \lim_{i \to \infty} \bs_i$ (since $\bs$ and $\bs_i$ share as prefix $\pal(v_i)$ whose length grows with $i$). 
By choice of $v_i$, we know that $\max\{ j \mid \alpha^j \in F(v_i^\omega), \alpha \in \A\} \leq \ell.$
By Theorem 5.2 in \cite{JustinPirillo2002} each $\bs_i$ is $(\ell+3)$-th power-free. Consequently $\bs$ is also $(\ell+3)$-th power-free.

Now by hypotheses, $\Delta = waua^\ell v a\Delta'$ with $a \in \A$ and $u, v \in \A^+$ such that $|u|_a=|v|_a=0$. Let $\bs'$ be the standard episturmian word directed by $va\Delta'$. The letter $a$ occurs in $\bs'$ and considering $b$ the first letter of $v$, we see that $b \neq a$ and $ab$ is a factor of the infinite word $\bs'$. 
Since $\ipal(\bw) = \lim_{n \to \infty} L_{a_1\ldots a_n}(a_{n+1})$ \cite[Cor. 2.7]{JustinPirillo2002},  $\bs$ contains as a factor the word $L_{waua^\ell}(ab) =$ $L_{wau}(a^{\ell+1}b)$ and so $\bs$ contains 
$L_{wa}(L_u(a)^{\ell+1}\pal(u)b)$. By  \cite{JustinPirillo2002}, since $a$ does not occur in $u$, $L_u(a) = \pal(u)a$ (since $a$ does not occur in $u$). Consequently
$L_{wa}(L_u(a)^{\ell+1}\pal(u)b) =$ $L_{wa}( (\pal(u)a)^{\ell+1}\pal(u)b)
=$\\ \hspace*{\fill} $L_w( L_a(\pal(u)a)^{\ell+1}$ $L_a(\pal(u))ab) =$ $L_w(L_{au}(a)^{\ell+2}b)$.\\
Hence $\bs$ contains the $(\ell+2)$-th power $(L_{wau}(a))^{\ell+2}$.
\end{proof}

Previous proposition can be viewed as a generalization of Lemma~\ref{L:power-free}. As a direct consequence, we have:

\begin{cor} \label{h:powers} The words $\bw_{\h,n}$ are $(n+4)$-th power free but contain $(n+3)$-th powers.
\end{cor}

Using previous results one can deduce the following properties of the words ${\bf s}_{\h,n,a,b,c}$.

\begin{proposition} Let ${\bf s}_{\h,n,a,b,c}$ be a fixed point of the $\ipal_\h$ operator, for a fixed $n$. Then ${\bf s}_{\h,n,a,b,c}$  satisfies the following properties:
\begin{enumerate}
\item It is not an episturmian word, but is a pseudostandard word.
\vspace{-0.2cm}
\item It is not a fixed point for some non trivial morphism.
\vspace{-0.2cm}
\item It is $(n+4)$-th power-free but contains $(n+3)$-th powers.
\vspace{-0.2cm}
\item The frequencies of the letters $b$ and $c$ are equal.
\vspace{-0.2cm}
\end{enumerate}
\end{proposition}

\section{Conclusion}

Let us summarize three problems raised by the content of this paper.

It is easy to see that any infinite word which is $k$-th power-free for an integer $k$ has a critical exponent. 
This is the case of all words studied in this paper. An open question is to find closed formulas of the values of the critical exponent of words ${\bf s}_{R,n,a,b}$, $\bs_{\h, n, a, n}$ and $\bs_{\e,a,b}$. Another direction of research would be to find a geometric interpretation of the palindromic closure. It may help find more about the fixed points of the operation we considered here. Finally since the study of the pseudostandard words which are fixed points of the $\pal_\V$ operator raises numerous intriguing questions, it might be interesting to also work with the more general families of words introduced in \cite{DeLuca2008}. The first one is called the {\it generalized pseudostandard words}, that is the pseudostandard words directed by two words: the traditional directive word and a word describing the antimorphism to used at each iteration. The second one is the pseudostandard words with seeds, that is the words obtained by iteration of the ${\oplus_\V}$ operator with a non empty word, called the seed.





%

{\footnotesize
\bibliographystyle{plain}
\bibliography{JPRV}

\begin{thebibliography}{10}

\bibitem{AdamczewskiBugeaud2007}
B.~Adamczewski and Y.~Bugeaud.
\newblock Palindromic continued fractions.
\newblock {\em Ann. Inst. Fourier (Grenoble)}, 57(1557-1574), 2007.

\bibitem{ArnouxRauzy1991}
P.~Arnoux and G.~Rauzy.
\newblock Repr\'esentation g\'eom\'etrique de suites de complexit\'e $2n+1$.
\newblock {\em Bull. Soc. Math. France}, 119(2):199--215, 1991.

\bibitem{CrispMoran1993}
D.~Crisp, W.~Moran, A.~Pollington, and P.~Shiue.
\newblock Substitution invariant cutting sequences.
\newblock {\em J. Th. Nombres Bordeaux}, 5:123--137, 1993.

\bibitem{deLuca1997}
A.~de~Luca.
\newblock Sturmian words: structure, combinatorics, and their arithmetics.
\newblock {\em Theor. Comput. Sci.}, 183:45--82, 1997.

\bibitem{DeLuca2008}
A.~De~Luca.
\newblock {\em Combinatorial aspects of {S}turmian sequences and their
  generalizations}.
\newblock Ph. d. thesis, Universit\`a degli Studi di Napoli Federico II, 2008.

\bibitem{deLucaDeLuca2006}
A.~de~Luca and A.~De~Luca.
\newblock Pseudopalindrome closure operators in free monoids.
\newblock {\em Theor. Comput. Sci.}, 362(1-3):45--82, 2006.

\bibitem{DroubayJustinPirillo2001}
X.~Droubay, J.~Justin, and G.~Pirillo.
\newblock Episturmian words and some constructions of {de Luca} and {Rauzy}.
\newblock {\em Theor. Comput. Sci.}, 255:539--553, 2001.

\bibitem{GlenJustin2009TIA}
A.~Glen and J.~Justin.
\newblock Episturmian words: a survey.
\newblock {\em RAIRO Theoretical Informatics and Applications}, To appear.

\bibitem{JustinPirillo2002}
J.~Justin and G.~Pirillo.
\newblock Episturmian words and episturmian morphisms.
\newblock {\em Theor. Comput. Sci.}, 276(1-2):281--313, 2002.

\bibitem{Kolakoski1965}
W.~Kolakoski.
\newblock Self {G}enerating {R}uns, {P}roblem 5304.
\newblock {\em Amer. Math. Monthly}, 72:674, 1965.

\bibitem{Lothaire2002}
M.~Lothaire.
\newblock {\em Algebraic Combinatorics on Words}, volume~90 of {\em
  Encyclopedia of Mathematics and its Applications}.
\newblock Cambridge University Press, 2002.

\bibitem{Shur2000}
A.M. Shur.
\newblock The structure of the set of cube-free z-words in a two-letter
  alphabet.
\newblock {\em IZV Math}, 64(4):847--861, 2000.

\bibitem{Vandeth2000}
D.~Vandeth.
\newblock Sturmian words and words with a critical exponent.
\newblock {\em Theor. Comput. Sci.}, 242:283--300, 2000.

\end{thebibliography}
}

\end{document}